\documentclass[a4paper,12pt]{amsart}
\usepackage{amssymb,amsthm,amsmath,color,url}
\usepackage[margin=3cm]{geometry}

\usepackage{graphicx}
\usepackage{stmaryrd}
\usepackage[T1]{fontenc}
\usepackage[english]{babel}
\usepackage{comment}
\usepackage{hyperref}
\usepackage{ifthen}
\hypersetup{
  colorlinks   = true,
  linkcolor = red!70!black,
     citecolor    = green!50!black
}
\usepackage[noabbrev,capitalise]{cleveref}
\crefname{subsection}{Subsection}{Subsections}
\crefname{claim}{Claim}{Claims}
\crefname{problem}{Problem}{Problems}

\makeatletter
\def\namedlabel#1#2{\begingroup
   \def\@currentlabel{#2}%
   \label{#1}\endgroup
}
\makeatother

\usepackage{thmtools}
\usepackage{thm-restate}
\usepackage{enumitem}
\usepackage{caption,subcaption}

\usepackage{tikz}
\usetikzlibrary{shadows}
\usetikzlibrary{backgrounds}
\usetikzlibrary{arrows}
\usetikzlibrary{patterns}
\usetikzlibrary{calc}
\usetikzlibrary{shapes}
\usetikzlibrary{positioning}
\usetikzlibrary{math}
\usetikzlibrary{external}
\usetikzlibrary{decorations.pathreplacing}
\declaretheorem[name=Theorem, numberwithin=section]{theorem}

\declaretheorem[name=Corollary, sibling=theorem]{corollary}

\declaretheorem[name=Claim, sibling=theorem]{claim}
\declaretheorem[name=Claim, numbered=no]{claim*}

\def\cqedsymbol{\ifmmode$\lrcorner$\else{\unskip\nobreak\hfil
\penalty50\hskip1em\null\nobreak\hfil$\lrcorner$
\parfillskip=0pt\finalhyphendemerits=0\endgraf}\fi}

\newcommand{\Cay}{\mathrm{Cay}}

\DeclareMathOperator{\girth}{girth}
\DeclareMathOperator{\diam}{diam}

\let\le\leqslant
\let\ge\geqslant

\let\geq\geqslant

\begin{document}
	
	\title[Optimization in graphical small
        cancellation theory]{Optimization in graphical\\
          small cancellation theory}

\author[L.~Esperet]{Louis Esperet}
\address[L.~Esperet]{Univ.\ Grenoble Alpes, CNRS, Laboratoire G-SCOP,
  Grenoble, France}
\email{louis.esperet@grenoble-inp.fr}

\author[U.~Giocanti]{Ugo Giocanti}
\address[U.~Giocanti]{Univ.\ Grenoble Alpes, CNRS, Laboratoire G-SCOP,
  Grenoble, France}
\email{ugo.giocanti@grenoble-inp.fr}

\thanks{The authors are partially supported by the French ANR Project
  GrR (ANR-18-CE40-0032), TWIN-WIDTH
  (ANR-21-CE48-0014-01), and by LabEx
  PERSYVAL-lab (ANR-11-LABX-0025).}

\begin{abstract}
Gromov (2003) constructed finitely generated groups whose
     Cayley graphs contain all graphs from a given infinite
     sequence of expander graphs of unbounded girth and bounded
     diameter-to-girth ratio. These so-called \emph{Gromov monster
       groups} provide examples of finitely generated groups that do
    not coarsely embed into Hilbert space, among other interesting
    properties. If graphs in Gromov's construction admit graphical small cancellation labellings, then one gets similar examples of Cayley graphs containing all the graphs of the family as isometric subgraphs. Osajda (2020) recently showed how to obtain such
    labellings using the probabilistic method. 
    In this short note, we simplify Osajda's
    approach, decreasing the number of generators of the
    resulting group significantly.
     \end{abstract}
        
	\maketitle

        \section{Introduction}

Given a finitely generated group $\Gamma$, with a finite set $S$ of
generators such that $S^{-1}=S$, the Cayley graph $\Cay(\Gamma,S)$ is the graph whose
vertices are the elements of $\Gamma$, in which we add an edge between
$\gamma$ and $\gamma \cdot s$ for any $\gamma\in \Gamma$ and $s\in
S$. Cayley graphs are a central object of study in geometric group
theory. It turns out that a number of
interesting properties of a group $\Gamma$ do not depend of the choice of the
generating set $S$. In particular, in order to show that $\Gamma$
\emph{does not satisfy} a given property of this type, it is
sufficient to find one generating set $S$ such that the corresponding
Cayley graph $\Cay(\Gamma,S)$ has a pathological behaviour.

Consider a sequence $\mathcal{G}=(G_n)_{n\ge 1}$ of bounded
degree graphs, whose girth (length of a shortest
non-trivial cycle) tends to infinity. We say that the
sequence is  \emph{dg-bounded} if
the ratio between the diameter and the girth of each $G_n$ is bounded by a
(uniform) constant, see \cite{AT18}. Consider such a sequence
$\mathcal{G}$. Gromov \cite{Gro03} proved that there is a 
finitely generated group $\Gamma$ whose Cayley graph contains (in a certain
metric sense) all the members of $\mathcal{G}$. By choosing
$\mathcal{G}$ as a family of suitable expander graphs, this implies that such a
group $\Gamma$ has a number of pathological properties, in particular
related to coarse embeddings in Hilbert space, or to Guoliang Yu's
property A. The construction has also been used very recently to disprove a
conjecture on the twin-width of groups and hereditary graph classes
\cite{BGTT22}. Gromov \cite{Gro03}  introduced the graphical small
cancellation condition on the labellings. By the classical small
cancellation theory, the existence of labellings of
$\mathcal{G}=(G_n)_{n\ge 1}$ with the graphical small cancellation
condition guarantees that in Gromov’s construction each graph $G_n$
embeds \emph{isometrically} in the Cayley graph $\Cay(\Gamma,S)$,  which means that the embedding
of each $G_n$ in  $\Cay(\Gamma,S)$ is distance-preserving and thus in
particular the
graphs $G_n$ appear as induced subgraphs in $\Cay(\Gamma,S)$. Osajda
\cite{Osajda} recently showed, using the probabilistic method, that
the graphical small cancellation labellings do exist, under mild
assumptions on $\mathcal{G}=(G_n)_{n\ge 1}$.


%
%
Given a sequence $\mathcal{G}=(G_n)_{n\ge 1}$ of graphs whose edges
are labelled with elements from some set $S$, a \emph{word} in 
$\mathcal{G}$ is a sequence of labels that can be read along a path of
some graph of $\mathcal{G}$.
The main idea of graphical small cancellation theory is to assign labels
from a finite set $S$ to the edges of all the graphs from the sequence
$\mathcal{G}=(G_n)_{n\ge 1}$, such that words in each
$G_n$ that are sufficiently long compared to the girth of $G_n$ occur
only once in all the sequence $\mathcal{G}$ (this will be made more
precise in the next section). 
The labels from $S$ are then used as generators to
define the group $\Gamma$ whose relators are the words labelling the
cycles of each $G_n$. The number of labels (the size of the set $S$)
then gives an upper bound to the minimum number of generators of the group, and thus the
degree of the associated Cayley graph (up to a multiplicative factor
of two, if we do not require that $S$ is closed under taking inverses). A natural problem  is to
minimize this number of generators.

The purpose of the present note is twofold: we present a
simplified version of the proof of existence of the labelling of  Osajda
\cite{Osajda}, and significantly decrease the number of generators
(and thus the degree of the corresponding Cayley graph). Osajda's
proof is based on an application of the Lov\'asz
  Local Lemma. Instead, we use a self-contained counting argument
  popularized by Rosenfeld \cite{Rosenfeld20}, and originally introduced
  in the field of combinatorics on words in the context of pattern avoidance. This allows us to
  cleanly handle all the different forbidden patterns at once, instead of sequentially, and
  greatly reduces the number of labels. We combine this with a
  significantly simpler (and stronger) analysis of intersecting patterns in order to a obtain a shorter argument
  that also produces much better bounds.

  For the sake of  concreteness, if we take $\mathcal{G}=(G_n)_{n\ge 1}$ to be the
  sequence of cubic Ramanujan graphs introduced by Chiu \cite{Chiu92}, which
  is likely to offer the best known parameters in terms of degree and
  diameter-to-girth ratio, our result leads to the existence of a
  group with 96 generators, whose Cayley graph (of maximum degree 96) contains all the
  graphs from  $\mathcal{G}$ as isometric subgraphs. For the same family, the construction of Osajda \cite{Osajda} uses about
  $10^{272}$ generators (although we note that some of the quick
  optimization steps we perform in Section \ref{sec:optim} can also be
  carried directly in Osajda's proof, improving his bound to about
  $10^{70}$ generators).

  \section{Preliminaries}

  All the graphs we consider in the paper are initially
  undirected. Each graph $G$ is then given  an arbitrary
  orientation $\vec{G}$ (i.e., the choice of a direction, for each edge of
  $G$). The results do not depend on
  the specific orientation, but the orientation is nevertheless
  crucial to define the relevant objects that we consider belows. Consider a
  set $S$ which is closed under (formal) inverse (that is, there is an
  involution without fixed point between the elements of $S$, which
  we denote by $a\mapsto \bar{a}$). Consider also a
  labelling $\ell: E(G)\to S$ of the edges of $G$ by the elements of
  $S$. We extend the labelling $\ell$ to the ordered pairs of adjacent
  vertices $(x,y)$ in $G$ as follows: if $(x,y)$ is an arc of
  $\vec{G}$ then $\ell(x,y)=\ell(xy)$ and otherwise
  $\ell(x,y)=\overline{\ell(xy)}$. The orientation $\vec{G}$ is only
  used to define this extended labelling $\ell$ of the the ordered pairs of adjacent
  vertices, and will not be mentioned elsewhere. 
  We say that the labelling $\ell$ is \emph{reduced} if for any vertex
  $v\in V(G)$, and for any pair of distinct neighbors $u,w$ of $v$ in
  $G$, $\ell(v,u)\ne \ell(v,w)$. An \emph{$\ell$-word} (or simply a
  \emph{word}, if $\ell$ is clear from the context) in $G$ is obtained from a path $P$
  in $G$ as follows: if $P=v_1,v_2,\ldots,v_k$, then
  $\ell(P):=\ell(v_1,v_2)\cdots \ell(v_{k-1},v_k)\in L^*$ is the
  $\ell$-word associated to $P$. The \emph{length} of a path is its
  number of edges. We remark that in this paper we consider paths as
  either a sequence of vertices, or a sequence of edges, depending on
  the context, and in particular any path $P=v_1,v_2,\ldots,v_k$ is
  distinct from the reverse path $\overleftarrow{P}:=v_k,v_{k-1},\ldots,v_1$.

  \medskip

  The girth (length of a smallest cycle) of a graph $G$ is denoted by
  $\girth(G)$, and its diameter (the maximum distance between
  two vertices of $G$) is denoted by $\diam(G)$.
  Let $\mathcal{G}=(G_n)_{n\ge 1}$ be a sequence of graphs. Let
  $\lambda$ be a positive real number (for the main application in
  group theory we need $\lambda\in \left(0,\tfrac16\right]$, but this
  will not be needed in the full generality of the results presented
  in this section and the next). Following the terminology of \cite{Osajda}, a sequence
  of labellings $(\ell_n)_{n\ge 1}$ of the graphs from $\mathcal{G}$,
  with labels from some set $S$ as above, is
  said to satisfy the \emph{$C'(\lambda)$-small cancellation property} if for
  all $n\ge 1$, $\ell_n$ is a reduced labelling of $G_n$ and no word of length at
  least $\lambda \cdot \girth(G_n)$ in $G_n$ appears on a different path in
  $\mathcal{G}$. 
  Small cancellation properties were initially introduced for groups,
  as a convenient tool to construct \emph{word-hyperbolic groups}, see
  for instance Chapter V in \cite{LS01}. The
  property $C'(\lambda)$ we use here is defined in the more general
  context of graphs, and is usually known as \emph{graphical
    cancellation property} in the literature. In the remainder of the
  paper we will omit the ``graphical'' term, as there is no risk
  of confusion with the original small cancellation properties.


  Osajda~\cite{Osajda} recently proved that under mild assumptions,
  any sequence of bounded degree dg-bounded graphs of unbounded girth admits small cancellation labellings with a
  finite number of labels.

  \begin{theorem}[\cite{Osajda}]\label{thm:Osajda}
    Let $\lambda\in \left(0,\tfrac16\right]$ and $A>0$ be real numbers, and let $\Delta\ge 3$ be an integer. 
    Let $\mathcal{G}=(G_n)_{n\ge 1}$ be a sequence of graphs of
    maximum degree $\Delta$ such that
    $\girth(G_n)\to \infty$ as $n\to \infty$, and
    $\diam(G_n)\le A \cdot \girth(G_n)$ for any $n\ge
    1$. Assume moreover that $1<\lfloor \lambda \cdot
    \girth(G_n)\rfloor <\lfloor \lambda \cdot 
    \girth(G_{n+1})\rfloor$ for every $n\ge 1$.     Let $$L\ge
    2e^4\Delta^{2A/\lambda+2}\cdot(4e^4\Delta)^{8A/\lambda+16}$$ be any
    even integer. Then $\mathcal{G}$ has a sequence of labellings satisfying the
    $C'(\lambda)$-small cancellation property, with labels from a set
    $S$ of size $L$.
  \end{theorem}


The bound on  $L$ in Theorem
  \ref{thm:Osajda} has two components: $2e^4\Delta^{2A/\lambda+2}$
  comes from a first phase, where Osajda shows how to assign labels in
  each $G_n\in \mathcal{G}$, so that no word of $G_n$ appears as a
  word of length at least $\lambda\cdot \girth(G_i)$ in some $G_i$,
  with $i<n$. The second component, $(4e^4\Delta)^{8A/\lambda+16}$,
  comes from a second phase where Osajda shows how to  assign labels in
  each $G_n\in \mathcal{G}$, so that no word of $G_n$ of length at
  least $\lambda\cdot \girth(G_n)$ appears twice in $G_n$. This second
  phase is significantly more involved, which explains the much larger
  label size. Our main contribution is the following.

  \begin{itemize}
    \item we use a counting argument instead of the Lov\'asz Local
      Lemma. This allows us to assign labels in a single phase
      (resulting in an additive combination of the number of labels, instead
      of a multiplicative one), and optimize the multiplicative
      constants. Moreover, the resulting proof is completely
      self-contained.
      \item we provide a major simplification in the analysis of
        Osajda's second phase, showing that the long words appearing twice in $G_n$ can
        be avoided with a number of labels of size comparable to
        Osajda's first phase.
      \end{itemize}
      
\medskip
  
  Using results of Gromov (see \cite{Ol03,Gru15}),  Theorem
  \ref{thm:Osajda} leads to the
  following.

  \begin{corollary}\label{cor:Osajda}
Let $\lambda,A,\Delta,\mathcal{G}$ be
as in Theorem \ref{thm:Osajda}. Then for any even integer $L\ge
2e^4\Delta^{2A/\lambda+2}\cdot(4e^4\Delta)^{8A/\lambda+16}$, 
there is a group $\Gamma$ with a set $S$ of $L$ generators such that the
corresponding Cayley graph $\Cay(\Gamma,S)$ contains isometric copies
of all the graphs from $\mathcal{G}$.
\end{corollary}

As alluded to in the introduction, in applications we typically want $\mathcal{G}$ to be a sequence of
\emph{expander graphs}. We omit the precise definition here, as it will not
be necessary in this paper. We only mention that expansion can be
defined in several essentially equivalent ways, using isoperimetric
inequalities or spectral properties. Families of random regular graphs
typically have these properties, but constructing explicit families of
expander graphs has been an important problem in Mathematics, with
major applications in Theoretical Computer Science. We refer the
interested reader to the survey \cite{HLW06} for more on expander
graphs.

A
useful family $\mathcal{G}$ for us is the sequence of cubic
Ramanujan graphs introduced by Chiu \cite{Chiu92}. These graphs are
expander graphs (as Ramanujan graphs, they have the best possible spectral
expansion), are $\Delta$-regular with $\Delta=3$, satisfy $\girth(G_n)\to\infty$
as $n\to \infty$ (their girth is logarithmic in their number
of vertices) and $\diam(G)\le
\tfrac32\, \girth(G)+5$ for any $G\in \mathcal{G}$. By discarding a
bounded number of small graphs in the sequence, this implies that we have $\diam(G)\le
(\tfrac32+\epsilon) \girth(G)$ for any $\epsilon>0$ and any graph $G$
in the sequence, and thus we can
take $A\le \tfrac32+\epsilon$ for any $\epsilon>0$.
  
  \section{Smaller cancellation labellings}

  Our main result is the following optimized version of Theorem
  \ref{thm:Osajda}.

  \begin{theorem}\label{thm:main}
    Let $\lambda,A,\Delta,\mathcal{G}$ be as in Theorem~\ref{thm:Osajda},
    that is $\lambda\in \left(0,\tfrac16\right]$ and $A>0$ are real numbers,  $\Delta\ge 3$ is an
    integer, and $\mathcal{G}=(G_n)_{n\ge 1}$ is a sequence of graphs of
    maximum degree $\Delta$ such that
    $\girth(G_n)\to \infty$ as $n\to \infty$, and
    $\diam(G_n)\le A \cdot \girth(G_n)$ and  $1<\lfloor \lambda \cdot
    \girth(G_n)\rfloor <\lfloor \lambda \cdot 
    \girth(G_{n+1})\rfloor$ for every $n\ge 1$.
    Let  $$L\ge
2(\Delta-1)+26 (\Delta-1)^{2A/\lambda+2}$$ be any even
integer. Then $\mathcal{G}$ has a sequence of labellings satisfying the
    $C'(\lambda)$-small cancellation property, with labels from a set
    $S$ of size  $L$.
  \end{theorem}

We note that the
multiplicative constant of 26 in the bound on $L$ can be optimized both for small values of
$\Delta$ and asymptotically as $\Delta\to \infty$. We have chosen not
to do so here
for simplicity, and we remark that improving the factor 2 in the
exponent of $(\Delta - 1)$ is a more rewarding challenge (see the next
section). When $\Delta\to \infty$, the number $L$ of labels in
Theorem \ref{thm:main} grows
as $O(\Delta^{2A/\lambda+2})$, and we will see in the next section
that this can be easily improved to $O(\Delta^{A/\lambda+2})$. This is
to be compared with the bound
$O(\Delta^{10A/\lambda+18})$ of Theorem
  \ref{thm:Osajda}. In the next section we will also see several
  ways to improve the constants significantly when $\Delta=3$,
  and the girth of the first graph in the sequence is already quite large.

\medskip

  Similarly as above, we obtain the following corollary.

  \begin{corollary}\label{cor:main}
Let $\lambda,A,\Delta,\mathcal{G},L$ be as in Theorem~\ref{thm:main}. Then
there is a group $\Gamma$ with a set $S$ of  $L$ generators such that the
corresponding Cayley graph $\Cay(\Gamma,S)$ contains isometric copies
of all the graphs from $\mathcal{G}$.
\end{corollary}

Using the family of cubic Ramanujan graphs of Chiu \cite{Chiu92} mentioned at the
end of the previous section, we can apply
Corollary \ref{cor:main} with $\Delta=3$, $A=\tfrac32$ and
$\lambda=\tfrac16$. Then we obtain a group with a set of 
$L=4+26\cdot 2^{20}=27262980$
generators such that the corresponding Cayley graph contains
isometric copies of graphs from an infinite family of expander
graphs. We will see in Section \ref{sec:optim} how to decrease this
number of generators to 96.

If instead we apply Corollary \ref{cor:Osajda} to the same family
$\mathcal{G}$ (and hence with the same parameters $\Delta=3$,
$A=\tfrac32$ and $\lambda=\tfrac16$), the
resulting Cayley graph has degree more than $10^{272}$.
\bigskip

We now prove our main result.

\bigskip

\noindent \emph{Proof of Theorem \ref{thm:main}.}
Let $\alpha:=2(\Delta-1)^{2A/\lambda+2}$, and
let  $$L\ge 2(\Delta-1)+13\alpha =2(\Delta-1)+26(\Delta-1)^{2A/\lambda+2} $$ be an
even integer.
Let $S$ be a set of $L$
elements, closed under formal inverses (and such that each element
$a\in S$ is
different from its formal inverse $\bar{a}$). For any $n\ge 1$,
let $\gamma_n:=\lfloor \lambda\cdot \girth(G_n)\rfloor$. In particular
$\gamma_n\le \lambda\cdot \girth(G_n)\le \gamma_n+1$ for any $n\ge 1$,
and thus

\begin{equation}\label{eq:1}
  \frac1\lambda\le \frac{\girth(G_n)}{\gamma_n}\le
  \frac1\lambda+\frac1{\lambda\gamma_n}\le \frac2\lambda.
\end{equation}

We will sequentially assign labels from $S$ to the edges of each of the graphs
$(G_n)_{n\ge 1}$. Assume that for each $i<n$, we have already defined a
labelling $\ell_i$ of the edges of
$G_i$  such that the sequence of labellings $(\ell_i)_{i<n}$ satisfies the $C'(\lambda)$-small
cancellation property. We now want to define a labelling $\ell_n$ of
$G_n$ so that the sequence
$(\ell_i)_{i\le n}$ of labellings of the graphs from $(G_i)_{i\le n}$ still satisfies the $C'(\lambda)$-small
cancellation property.

For the proof it will be convenient to consider \emph{partial} labellings of
$G_n$, which are labellings of some subset $F$ of edges of
$G_n$. Equivalently, these are labellings of the edges of $G_n[F]$,
the subgraph of $G_n$ induced by the edges of $F$. We recall that each
labelling $\ell(xy)$ of an edge $xy$ yields two labellings $\ell(x,y)$
and $\ell(y,x)$ of the pairs $(x,y)$ and $(y,x)$ by elements of $S$
that are formal inverse (and that whether $\ell(xy)=\ell(x,y)$ or
$\ell(xy)=\ell(y,x)$ depends only on the orientation of the edge $xy$
in some fixed but otherwise arbitrary orientation of the graph under consideration).

\smallskip

Let $F$ be a non-empty subset of $E(G_n)$. We say that a
labelling $\ell$ of $G_n[F]$ with labels from $S$ is \emph{valid} if it satisfies the following
properties:
\begin{enumerate}
  \item[(a)] $\ell$ is a reduced labelling of $G_n[F]$,
    \item[(b)] for each $1\le i < n$, no $\ell_i$-word of length at least
      $\gamma_i$ in $G_i$  appears as an $\ell$-word in
      $G_n[F]$, and
      \item[(c)] no $\ell$-word of length at least $\gamma_n$ appears on two
        different paths of $G_n[F]$.
      \end{enumerate}

      Let  $c(F)$ be the number of valid
labellings $\ell$ of $G_n[F]$ with labels from $S$ (when $F$ is empty we conveniently define $c(F):=1$). 
  In the remainder of the proof we will show the following claim,
  which clearly implies that $G_n$ has a labelling $\ell_n$ such that
  the sequence of labellings $(\ell_i)_{i\le n}$ of $(G_i)_{i\le n}$ still satisfies the $C'(\lambda)$-small
cancellation property, and thus we can find such labellings in all
the graphs from $\mathcal{G}$.

      \begin{claim}\label{cl}
For any non-empty $F\subseteq E(G_n)$ and any $e\in F$, $c(F)\ge
\alpha\cdot c(F\setminus\{e\})$.
\end{claim}

We prove the claim by induction on $|F|$. Recall that by assumption,
$\gamma_i>1$ for any $i\ge 1$, so the properties (a), (b), (c) above
are trivially satisfied if $F$ contains a single element $e$, which is
assigned an arbitrary label from $S$. It
follows that $c(\{e\})=L\ge \alpha=\alpha\cdot c(\emptyset)$, as
desired. So we can now assume that $F$ contains at least two
elements.

Assume that we have proved the claim for any $F'\subseteq E(G_n)$
with $|F'|<|F|$. Consider any edge $xy\in F$. Our goal in the
remainder of the proof is to show
that $c(F)\ge
\alpha\cdot c(F\setminus\{xy\})$. Note that by the induction hypothesis, for any
subset $F'\subseteq F$ containing $xy$, 
\begin{equation}\label{eq:ind}
  c(F\setminus F')\le
  \alpha^{1-|F'|}\cdot c(F\setminus \{xy\}).
\end{equation}
Let $\mathcal{L}$ denote
the set of  labellings $\ell$ of $F$ with labels from $S$
whose restriction to $F\setminus \{xy\}$ is valid, but such that
$\ell$ itself is not. Then
\begin{equation}\label{eq:3}
  c(F)=L\cdot c(F\setminus\{xy\})-|\mathcal{L}|.
\end{equation}
Consider first the subset $\mathcal{L}_a\subseteq \mathcal{L}$ of labellings of $F$ that do not satisfy  (a) above. Then by definition, for
any $\ell\in \mathcal{L}_a$, $x$
has a neighbor $z$ different from $y$ such that $\ell(x,y)=\ell(x,z)$,
or $y$
has a neighbor $z$ different from $x$ such that
$\ell(y,x)=\ell(y,z)$. By assumption, the labelling $\ell^-$ of $F\setminus\{xy\}$ obtained
from $\ell$ by discarding the label of $xy$ is valid. Moreover, $\ell$
can be recovered in a unique way from $\ell^-$ and the edge $xz$ or $yz$ as above. As
there are at most $2(\Delta-1)$ choices for such an edge incident to
$xy$, we obtain
\begin{equation}\label{eq:la}
|\mathcal{L}_a|\le 2(\Delta-1)\cdot c(F\setminus\{xy\}).
\end{equation}
For $1\le i \le n-1$, let $\mathcal{L}_i$ be the subset of labellings
$\ell\in \mathcal{L}$ of $F$ such that $G_n[F]$ contains a path $P$ containing $xy$
such that
$\ell(P)$ coincides with some $\ell_i$-word $\ell_i(Q)$ of length $\gamma_i$ in
$G_i$. Let $\mathcal{L}_n$ be the subset of labellings $\ell\in
\mathcal{L}\setminus \mathcal{L}_a$ of $F$ such that $G_n[F]$ contains a path $P$ containing $xy$
such that
$\ell(P)$ coincides with some $\ell$-word $\ell(Q)$ of length $\gamma_n$ in
$G_n[F]$, for some path $Q$ distinct from $P$.

For each $1\le i \le n$ and each
labelling
$\ell\in \mathcal{L}_i$ as above, let $\ell^-$ denote the labelling of
$F\setminus E(P)$ obtained from $\ell$ by discarding the labels of the
edges of $P$. Then $\ell^-$ is a valid labelling of $F\setminus
E(P)$. Moreover, if $1\le i \le n-1$ or if $i=n$ and $P$ and $Q$ are
disjoint, then  $\ell^-$ together with the paths $P$ in $G_n$ and $Q$
in $G_i$ (where each path is viewed as a sequence of edges) are
sufficient to recover $\ell$ in a unique way.

Assume now that $\ell\in \mathcal{L}_n$ (so in particular $\ell$ is
reduced), and the distinct paths $P$ and $Q$ of length
$\gamma_n$ in  $G_n[F]$ such that $\ell(P)=\ell(Q)$ are not
edge-disjoint. We first observe that $E(P)\cap E(Q)$ is a subpath of
$P$ and $Q$, since otherwise $G_n$ would contain a cycle of length
less than  $2\gamma_n$, contradicting the assumption that
$\girth(G_n)\ge \tfrac{\gamma_n}{\lambda}\ge 6 \gamma_n$. Let
$P=x_0,x_1,\ldots,x_{\gamma_n}$ and
$Q=y_0,y_1,\ldots,y_{\gamma_n}$. Then
$\ell(x_i,x_{i+1})=\ell(y_i,y_{i+1})$ for any $0\le i \le
\gamma_n-1$. Our goal is to show that despite the fact that the edges of
$E(P)\cap E(Q)$ have been unlabelled in $\ell^-$, we can still recover
$\ell$ from $\ell^-$, $P$ and $Q$.

\begin{figure}[htb]
 \centering
 \includegraphics{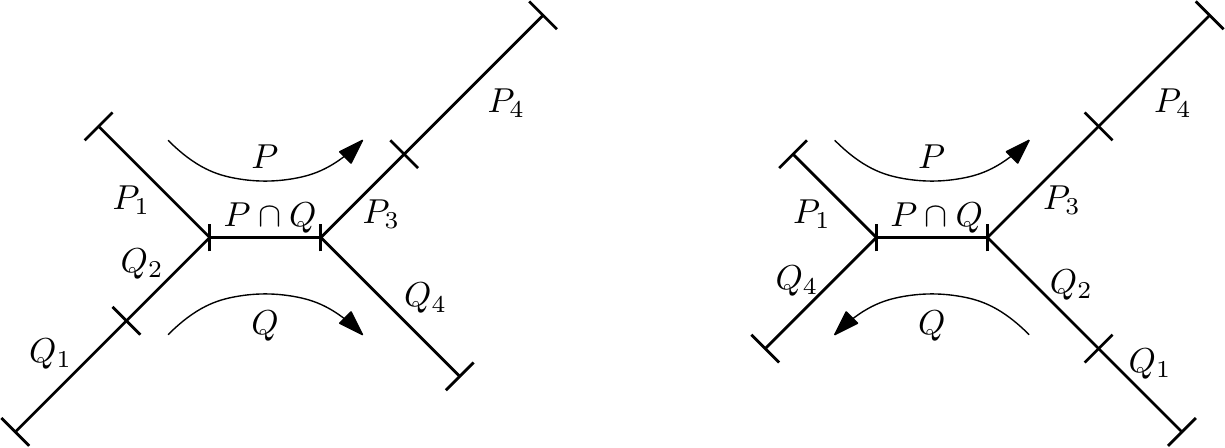}
 \caption{Two intersecting paths $P$ and $Q$.}
 \label{fig:paths}
\end{figure}

Assume first that $P$ and $Q$ intersect in the same direction, that is there are integers $0\le p,q\le
\gamma_n-1$ and $1\le k\le \gamma_n-1$ such that
$x_{p+i}=y_{q+i}$ for any $0\le i \le k$. Note that $p\ne q$ since
otherwise we would have $x_p=y_p$ and $x_{p+k}=y_{p+k}$ and the fact that $\ell(x_{p-1},x_p)=\ell(y_{p-1},y_p)$ or
$\ell(x_{p+k},x_{p+k+1})=\ell(y_{p+k},y_{p+k+1})$ would contradict the
fact that $\ell$ is reduced. Up to considering the reverse paths
$\overleftarrow{P}$ and $\overleftarrow{Q}$ instead of $P$ and $Q$, we
can assume without loss of generality that $q>p$. Divide  $P$ into
consecutive subpaths $P_1$, $P\cap Q$,  $P_3$, and $P_4$ and divide
$Q$ into consecutive subpaths $Q_1$, $Q_2$, $P\cap Q$, and $Q_4$, in such a way that
$\ell(P_1)=\ell(Q_1)$ and $\ell(P_4)=\ell(Q_4)$ (see Figure
\ref{fig:paths}, left). As $P_1$ and $P_4$ are edge-disjoint from $E(P)\cap E(Q)$, both
$\ell(P_1)$ and $\ell(P_4)$  can be recovered from $\ell^-$. Note that as we assumed that $q>p$, $|Q_2|>0$, i.e. $Q_2$ has at least one edge.
Let $P'$
be the subpath of $P$ obtained by concatenating $P\cap Q$ and $P_3$. It
remains to explain how to recover $\ell(P')$ from
$\ell^-$. For this, it suffices to observe that since
$\ell(P)=\ell(Q)$, the prefix of $\ell(P')$ of
size $|Q_2|$ must be equal to $\ell(Q_2)$. Then the prefix of $\ell(P')$ of
size $2|Q_2|$ must be equal to $\ell(Q_2) \cdot \ell(Q_2)$. 
By iterating
this observation, it follows
that $\ell(P')$ is a
prefix of the word $\ell(Q_2)^\omega$ (the concatenation of an
infinite number of copies of $\ell(Q_2)$). Since $Q_2$ is
edge-disjoint from $E(P)\cap E(Q)$, $\ell(P')$ (and thus $\ell(P)$) can be recovered from
$\ell^-$, $P$ and $Q$, as desired.

We now assume that $P$ and $Q$ intersect in reverse directions, that is there are integers $0\le p,q\le
\gamma_n$ and $k\ge 1$ such that
$x_{p+i}=y_{q-i}$ for any $0\le i \le k$. We say that $P$ and $Q$
\emph{collide} if there is an index $i$ such that either $x_i=y_i$, or $x_i=y_{i+1}$ and
$y_i=x_{i+1}$ (think of two particles following the trajectories of $P$ and $Q$ at the same
speed). Assume for the sake of contradiction that $P$ and $Q$ collide.
If $x_i=y_i$ for some index $i$, then $\ell$ is not reduced, which is
a contradiction. Otherwise we have
$\ell(x_i,x_{i+1})=\ell(y_{i},y_{i+1})=\ell(x_{i+1},x_i)$, which
contradicts the fact that
$\ell(x_i,x_{i+1})=\overline{\ell(x_{i+1},x_{i})}$ as for each $a\in S$, $\overline a \neq a$. So $P$ and $Q$ do
not collide, and in particular $p\ne q$. We recall that  $\overleftarrow{P}$ and  $\overleftarrow{Q}$ denote the paths
obtained by reversing $P$ and $Q$, respectively. When we use
this notation below we also write $\overrightarrow{P}$ and
$\overrightarrow{Q}$ instead of $P$ and $Q$ to avoid any confusion. Up
to considering $\overleftarrow{P}$ and $\overleftarrow{Q}$ instead of
$\overrightarrow{P}$ and $\overrightarrow{Q}$,  we can again assume without loss of generality that $q>p$. We divide $P$ into
consecutive subpaths $P_1$, $\overrightarrow{P}\cap \overleftarrow{Q}$, $P_3$ and
$P_4$ and we divide $Q$ into consecutive subpaths $Q_1$, $Q_2$,
$\overleftarrow{P}\cap \overrightarrow{Q}$, and $Q_4$, in such a way that
$\ell(P_1)=\ell(Q_1)$, $\ell(P_4)=\ell(Q_4)$ (see Figure \ref{fig:paths}, right). As
before,  $P_1$ and $P_4$ are edge-disjoint from $E(P)\cap E(Q)$, so both
$\ell(P_1)$ and $\ell(P_4)$  can be recovered from $\ell^-$. As $P$
and $Q$ do not collide, $|Q_2|>|\overrightarrow{P}\cap
\overleftarrow{Q}|$, which implies that $\ell(\overrightarrow{P}\cap
\overleftarrow{Q})$ is equal to a prefix of $\ell(Q_2)$, and can thus
be recovered from $\ell^-$. Finally, since $\ell(P)=\ell(Q)$,
$\ell(P_3)$ is equal to $\ell(\overleftarrow{P}\cap
\overrightarrow{Q})$, which is obtained by reading $\ell(\overrightarrow{P}\cap
\overleftarrow{Q})$ backwards. Hence, $\ell(P)$ can be recovered from
$\ell^-$, $P$ and $Q$, as desired.

\smallskip

For each $1\le i \le n$ and each edge $e$ in
$G_i$ there
are at most $(\Delta-1)^{\gamma_i-1}$ paths of length
$\gamma_i$ containing $e$ in which $e$ is at a fixed position on the
path. Hence, there are at at most $2\gamma_i(\Delta-1)^{\gamma_i-1}$ paths of length
$\gamma_i$ containing $e$ (and in particular at most
$2\gamma_i(\Delta-1)^{\gamma_i-1}$ choices for the path $P$ in $G_n$ containing $xy$ when
considering a labelling $\ell\in\mathcal{L}_i$).  Moreover, each $G_i$
has at most $1+\Delta+\Delta(\Delta-1)+\cdots
+\Delta(\Delta-1)^{\diam(G_i)-1}$ vertices, and thus at most
\begin{equation}\label{eq:edges}
  \tfrac{\Delta}2\cdot
  \left(1+\Delta\tfrac{(\Delta-1)^{\diam(G_i)}-1}{\Delta-2}\right)\le \tfrac32(\Delta-1)^{\diam(G_i)+2}
\end{equation}
edges, using $\Delta\ge 3$ (the inequality is quite loose here, we
have chosen the right-hand side mostly in order to simplify the computation later). It follows that each $G_i$ has at most
\begin{equation}\label{eq:edges2}
\tfrac32(\Delta-1)^{\diam(G_i)+2}\cdot 2(\Delta-1)^{\gamma_i-1}\le 3(\Delta-1)^{(2A/\lambda+1)\gamma_i+1}
\end{equation}
paths of length $\gamma_i$ (here the multiplicative factor $\gamma_i$
disappears since we can count each path from its starting edge). It
follows that there are at most $3
(\Delta-1)^{(2A/\lambda+1)\gamma_i+1}$ choices for the path $Q$ in
$G_i$ when
considering a labelling $\ell\in\mathcal{L}_i$.
Since $|E(P)|=\gamma_i$, it follows from \eqref{eq:ind} that for each labelling $\ell\in\mathcal{L}_i$, the number
of valid labellings $\ell^-$ of  $F\setminus E(P)$ is 
$c(F\setminus E(P))\le \alpha^{1-\gamma_i} \cdot c(F\setminus
\{xy\})$. As each $\ell\in\mathcal{L}_i$ can be recovered from $\ell^-$,
$P$ and $Q$ in a unique way, we obtain
\begin{eqnarray*}
  |\mathcal{L}_i|& \le & 2\gamma_i(\Delta-1)^{\gamma_i-1}\cdot 3
(\Delta-1)^{(2A/\lambda+1)\gamma_i+1}\cdot  \alpha^{1-\gamma_i} \cdot c(F\setminus
                       \{xy\})\\
& \le & 6 \gamma_i                        (\Delta-1)^{(2A/\lambda+2)\gamma_i}
                       \cdot\alpha^{1-\gamma_i} \cdot c(F\setminus
                       \{xy\})\\
  & \le & 6 \gamma_i (\alpha/2)^{\gamma_i} \cdot\alpha^{1-\gamma_i} \cdot c(F\setminus
          \{xy\})\\
  & \le & 6\alpha \cdot \gamma_i  (1/2)^{\gamma_i} \cdot c(F\setminus
                       \{xy\}),
\end{eqnarray*}
where we have used $\alpha=2(\Delta-1)^{2A/\lambda+2}$ in the third inequality.
As a consequence
\begin{equation}\label{eq:li}
  \sum_{i=1}^{n} |\mathcal{L}_i|\le 6 \alpha \sum_{i=1}^n\gamma_i  (1/2)^{\gamma_i}  c(F\setminus
                       \{xy\}) \le 12\alpha \cdot c(F\setminus
                       \{xy\}),
                     \end{equation}
                     where we have used $\sum_{j=1}^\infty j (1/2)^j=2$.
As $\mathcal{L}=\mathcal{L}_a\cup \bigcup_{i=1}^n \mathcal{L}_i$, it
follows from \eqref{eq:la} and \eqref{eq:li} that
\begin{eqnarray*}
  |\mathcal{L}|&\le& c(F\setminus \{xy\})\cdot
                     \left(2(\Delta-1)+12\alpha
                     \right)\\
  &\le & c(F\setminus \{xy\}) (L-\alpha),
\end{eqnarray*}
by the definition of $L$. By \eqref{eq:3}, we have 
\begin{eqnarray*}
  c(F)&=&L\cdot c(F\setminus \{xy\})-|\mathcal{L}|\\
      &\ge & L\cdot c(F\setminus \{xy\})-(L-\alpha)  c(F\setminus \{xy\})\\
             &\ge& \alpha \cdot c(F\setminus \{xy\}),
\end{eqnarray*}
as desired. This completes the proof of Claim \ref{cl}, which concludes the
proof of Theorem \ref{thm:main}.
\hfill $\Box$

\section{Optimizing the number of generators}\label{sec:optim}

So far our goal was to optimize the construction of Osajda \cite{Osajda},
while obtaining a result that is comparable to his (i.e., a result with
the exact same set of initial assumptions).
There are two quick ways to further optimize the number of labels in Theorem
\ref{thm:main}, if we have some control over the family
$\mathcal{G}$.

The first way consists in removing all sufficiently
small graphs from $\mathcal{G}$ (we have done this already with the
cubic Ramanujan graphs of Chiu \cite{Chiu92}, to argue that $A$ was
arbitrarily close to $\tfrac32$ in this case). As the girth of the graphs in
$\mathcal{G}$ tends to infinity, the right-hand-side of \eqref{eq:1}
can be replaced by $\tfrac{1+\epsilon}\lambda$ for any
$\epsilon>0$. This allows to replace all instances of $2A/\lambda$ by
$(1+\epsilon)A/\lambda$ in the proof, effectively
dividing by 2 the exponent of the number of labels in the
theorem. Using this observation in the case of the cubic Ramanujan
graphs of Chiu \cite{Chiu92}, with $\lambda=1/6$, we obtain $\alpha=2\cdot 
2^{(1+\epsilon)\tfrac32/\tfrac16+2}\le 4097$ for sufficiently small $\epsilon>0$, and a number of
labels $L\ge 2\cdot 2+13\cdot 4097\approx
53266$ is sufficient. 

\medskip

A more efficient way to decrease the number of labels in the case of
families of expander graphs with an explicit description consists in using a more precise bound on the number of
edges in a graph $G_n\in \mathcal{G}$, as a function of
$\girth(G_n)$. In \eqref{eq:edges}, we have used that $|E(G_n)|\le
\tfrac32(\Delta-1)^{\diam(G_n)+2}\le \tfrac32(\Delta-1)^{A\girth(G_n)+2} $. However, better bounds are known for a number of
    families $\mathcal{G}$. This is
    the case for the cubic Ramanujan graphs of Chiu \cite{Chiu92}
    mentioned in the previous section. The graphs $G$ in this class
    satisfy $|E(G_n)|\le \tfrac32\cdot 2^{(3 \girth(G_n)+6)/{4}}$, which is an
    improvement over the bound based on the diameter (recall that for
    these graphs $\Delta=3$ and $A$ can be made arbitrarily close to
    $\tfrac32$). Fix any real $\epsilon>0$, and recall that
    $\gamma_n=\lfloor \lambda\cdot \girth(G_n)\rfloor$.  Using as in the previous paragraph the fact that the girth of
    the graphs from $\mathcal{G}$ can be made arbitrarily large by
    discarding a constant number of graphs from the family, we can
    assume that $\gamma_n \epsilon>\gamma_1\epsilon$ is larger than
    any fixed constant, and thus
    $(3 \girth(G_n)+6)/{4}\le \tfrac{3+\epsilon}{4\lambda}\gamma_n$ and
    $|E(G_n)|\le \tfrac32\cdot 2^{\gamma_n\cdot (3+\epsilon)/4\lambda}$,  for
    any $n\ge 1$. With $\lambda=1/6$, we obtain $|E(G_n)|\le
    \tfrac32\cdot 2^{(9+\epsilon)\gamma_n/2}$, for any $n\ge 1$. Substituting this bound in \eqref{eq:edges2},
    we obtain that there are at most $ 3\cdot
    2^{(11+\epsilon)\gamma_n/2-1}$, paths of length $\gamma_n$ in
    $G_n$. Substituting this bound in the proof of Theorem
    \ref{thm:main}, and defining $\alpha:=(1+\epsilon)2^{(13+\epsilon)/2}$, we obtain the following.
\begin{eqnarray*}
|\mathcal{L}_i|& \le &2\gamma_i 2^{\gamma_i-1}\cdot 3\cdot 2^{(11+\epsilon)\gamma_i/2-1}
                       \cdot\alpha^{1-\gamma_i} \cdot c(F\setminus
                       \{xy\})\\
  & \le & \tfrac{3\alpha}2 \gamma_i \cdot 2^{(13+\epsilon)\gamma_i/2}\cdot \alpha^{-\gamma_i} \cdot c(F\setminus
                       \{xy\})\\
  & \le & \tfrac{3\alpha}2 \gamma_i\cdot (\tfrac1{1+\epsilon})^{\gamma_i} \cdot c(F\setminus
                       \{xy\}),
\end{eqnarray*}
As $\sum_{j=1}^\infty j \cdot
(\tfrac1{1+\epsilon})^{j}$ converges,
we can choose again $\gamma_1$
sufficiently large so that the truncated sum $\sum_{j=\gamma_1}^\infty j \cdot
(\tfrac1{1+\epsilon})^{j}$ is arbitrarily small (say smaller than
$\epsilon/(\tfrac{3\alpha}2)$). We obtain $  \sum_{i=1}^{n} |\mathcal{L}_i|\le \epsilon\cdot c(F\setminus
\{xy\})$, and the same computation as in the proof of Claim \ref{cl}
shows that any even number $L\ge 2\cdot
2+\epsilon+\alpha=\alpha+\epsilon+4$  of labels is
sufficient. Using $\alpha=(1+\epsilon)2^{(13+\epsilon)/2}$, and taking 
    $\epsilon>0$ sufficiently small, we can obtain that $L= 96$ labels
    are sufficient. 

    So, we obtain a group with a set $S$ of 96 generators
    whose Cayley graph $\Cay(\Gamma,S)$ contains infinitely many
    graphs of the sequence of cubic Ramanujan graphs as isometric subgraphs.

    \section{Conclusion}

The number $L$ of labels in Theorem
\ref{thm:main} is of order $O(\Delta^{2A/\lambda+2})$, as $\Delta\to
\infty$, and the remarks in the previous section improve this bound to
$O(\Delta^{A/\lambda+2})$. In typical applications, $A$ is a small
constant and the bound becomes $\Delta^{O(1/\lambda)}$. We now observe
that this is the right order of magnitude. If $G$ is a
$\Delta$-regular 
graph of girth $g$, then the ball of radius $g/2$ centered in any
vertex induces a tree, and thus for any $\lambda <1/2$, $G$ contains
$\Omega(\Delta^{g/2+\lambda g-1})$ paths of length $\lambda g$ (the ball
of radius $g/2$ centered in a vertex  contains
$\Omega(\Delta^{g/2})$ edges and each of them is the starting
point of $\Omega(\Delta^{\lambda g-1})$ paths of length $\lambda
g$). By the $C'(\lambda)$-small cancellation property, all these paths
must correspond to different words. As there are at most $L^{\lambda
  g}$ possible words of length $\lambda g$, we obtain $L^{\lambda
  g}=\Omega(\Delta^{g/2+\lambda g-1})$, and thus $L=\Omega(\Delta^{1/2\lambda+1-1/g})$.  As the girth of the graphs in our family is
unbounded, it follows that $L=\Omega(\Delta^{1/2\lambda+1})$,
which shows that the bound in Theorem \ref{thm:main} is fairly close to the
optimum (up to a small multiplicative factor in the exponent). It
remains an interesting problem to close the gap between the upper and
lower bounds, both in the case of small degree ($\Delta=3$) and
asymptotically as $\Delta\to \infty$.

\medskip

  It might also be interesting to consider other cancellation
  properties. For an integer $k\ge 1$, a family of labellings
  $(\ell_n)_{n\ge 1}$ of a graph family $\mathcal{G}=(G_n)_{n\ge 1}$
  satisfies the  \emph{$C(k+1)$-small cancellation
  property} if for any $n\ge 1$, $\ell_n$ is reduced and no cycle $C$ in $G_n$ can be divided
into $k$ paths $P_1, \ldots,P_k$ such that for each $1\le i \le k$,
the $\ell_n$-word associated to $P_i$ appears on a different
path in  $\mathcal{G}$. This condition is weaker than the $C'(1/k)$-small cancellation
  property, but nevertheless allows to construct finitely generated
  groups with interesting properties when $k\geq 7$ \cite{Gru15}. A natural problem
  is to obtain a version of Theorem \ref{thm:main} for $C(k)$-small
  cancellation, with an improved exponent.

  \medskip

  We conclude with some algorithmic remarks. Using the constructive proof
  of the Lov\'asz Local Lemma  by Moser and Tardos \cite{MT10}, the
  original proof of existence of the labelling given by Osajda
  \cite{Osajda} can be turned into an efficient algorithm computing
  the labels, by which we mean a randomized algorithm, running in
  polynomial time (in the size of $G_n$), and computing a $C'(\lambda)$-small
  cancellation labelling for the sequence of graphs $(G_i)_{1\le i
    \le n}$. As our main goal was to obtain a simple,
  self-contained proof of the existence of the labels, we chose to use
  counting rather than constructive techniques such as the entropy
  compression method (see \cite{GMP20}). It turns out that  our result can also be
  obtained with this type of techniques, at the cost of a longer and
  more technical analysis. 

  \subsection*{Acknowledgements}

  We thank Goulnara Arzhantseva for her comments on a previous version
  of this manuscript.
                
\bibliographystyle{plain}
\bibliography{biblio}

\end{document}